\documentclass{article}
\usepackage{amssymb,latexsym}
\usepackage{amsmath}
\usepackage{graphics}
\usepackage{graphicx}
\newtheorem{theorem}{Theorem}

\newtheorem{definition}[theorem]{Definition}

\newenvironment{proof}[1][Proof]{\textbf{#1.} }{\ \rule{0.5em}{0.5em}}
\numberwithin{equation}{section} \numberwithin{theorem}{section}
\numberwithin{corollary}{section} \numberwithin{lemma}{section}
\numberwithin{remark}{section} \numberwithin{notation}{section}

\begin{document}
\title{Commutative Quaternion Matrices}
\author{Hidayet H\"{u}da K\"{O}SAL and Murat TOSUN}

\maketitle

\begin{abstract}
In this study, we introduce the concept of commutative quaternions and commutative quaternion matrices. Firstly, we give some properties of commutative quaternions
and their Hamilton matrices. After that we investigate commutative quaternion matrices using properties of complex matrices. Then we define the complex adjoint
matrix of commutative quaternion matrices and give some of their properties.

\textbf{Mathematics Subject Classification (2000)}:  11R52;15A33

\textbf{Keywords}: Commutative quaternions, Hamilton matrices, commutative quaternion matrices.\\
\end{abstract}

\section{INTRODUCTION}

\noindent Sir W.R. Hamilton introduced the set of quaternions in 1853 \cite{1}, which can be represented as

\[Q = \left\{ {q = {t} + {x}i + {y}j + {z}k:\,\,{t},{x},{y},{z} \in R} \right\}.\]

\noindent Quaternion bases satisfy the following identities known as the

\[{i^2} = {j^2} = {k^2} =  - 1,\,\,\,\,ij =  - ji = k,\,\,\,jk =  - kj = i,\,\,\,ki =  - ik = j.\]

\noindent From these rules it follows immediately that multiplication of quaternions is not commutative. So, it is not easy to study quaternion algebra problems. Similarly,
it is well known that the main obstacle in study of quaternion matrices, dating back to 1936 \cite{2}, is the non-commutative multiplication of quaternions. In this
sense, there exists two kinds of its eigenvalue, called right eigenvalue and left eigenvalue.  The right eigenvalues and left eigenvalues of quaternion matrices
must be not the same.  The set of right eigenvalues is always nonempty, and right eigenvalues are well studied in literature \cite{3}.  On the contrary, left
eigenvalues are less known.  Zhang \cite{4}, commented that the set of left eigenvalues is not easy to handle, and few result have been obtained.

\noindent After the discovery of quaternions by Hamilton, Segre proposed modified quaternions so that commutative property in multiplication is possible \cite{5}. This number
system, called commutative quaternion. Commutative quaternions are decomposable into two complex variables \cite{6}. The set of commutative quaternions is $4 - $
dimensional like the set of quaternions. But this set contains zero-divisor and isotropic elements. As has been noticed, there is no theory of commutative
quaternion matrices.  However, commutative quaternion algebra theory is becoming more and more important in recent years and has many important applications to
areas of mathematics and physics \cite{7}-\cite{8}-\cite{9}.  In this paper, we study properties of commutative quaternion matrices for its applications.

\section{COMMUTATIVE QUATERNION ALGEBRA}

\noindent A set of commutative quaternions are denoted by

\[H = \left\{ {q = t + ix + jy + kz:\,\,t,x,y,z \in R,\,\,i,j,k \notin R} \right\}\]

\noindent where the bases elements  $i,j,k$ satisfy the following multiplication rules:

\[{i^2} = {k^2} =  - 1,\,\,\,\,{j^2} = 1,\,\,\,\,ijk =  - 1,\,\,\,ij = ji = k,\,\,jk = kj = i,\,\,\,ki = ik =  - j.\]

\noindent It is obvious that the multiplication rules in set of commutative quaternions are commutative i.e. $xy = yx$  if   and  are commutative quaternions.

\noindent Multiplication of any commutative quaternions $q =  t + ix + jy + kz$  and $ {q_1} = {t_1} + i{x_1} + j{y_1} + k{z_1}$  are defined in the following ways,

\begin{equation}\label{1}
\begin{array}{l}
q{q_1} = t{t_1} - x{x_1} + y{y_1} - z{z_1} + i\left( {x{t_1} + t{x_1} + z{y_1} + y{z_1}} \right)\\
\,\,\,\,\,\,\, + j\left( {t{y_1} + y{t_1} - x{z_1} - z{x_1}} \right) + k\left( {z{t_1} + t{z_1} + x{y_1} + y{x_1}} \right)
\end{array}
\end{equation}

\noindent Equation (\ref{1}) can be represented by means of matrix multiplication. The representation as a $4 \times 4$  real matrix is

\begin{equation}\label{2}
\left( {\begin{array}{*{20}{c}}
t&{ - x}&y&{ - z}\\
x&{\,\,\,\,t}&z&{\,\,\,\,y}\\
y&{ - z}&t&{ - x}\\
z&{\,\,y}&x&{\,\,t}
\end{array}} \right)
\end{equation}

\noindent which is a useful way to compute quaternion

\[\left( {\begin{array}{*{20}{c}}
{{r_0}}\\
{{r_1}}\\
{{r_2}}\\
{{r_3}}
\end{array}} \right) = \left( {\begin{array}{*{20}{c}}
t&{ - x}&y&{ - z}\\
x&{\,\,\,t}&z&{\,\,\,y}\\
y&{ - z}&t&{ - x}\\
z&{\,\,y}&x&{\,\,\,t}
\end{array}} \right)\left( {\begin{array}{*{20}{c}}
{{t_1}}\\
{{x_1}}\\
{{y_1}}\\
{{z_1}}
\end{array}} \right)\]

\noindent where $q{q_1} = {r_0} + i{r_1} + j{r_2} + k{r_3}.$

\noindent For  $q = t + ix + jy + kz \in H,$ there exists three kinds of its conjugate, called principal conjugations:

\[\begin{array}{l}
{q^{\left( {\rm{1}} \right)}} = t - ix + jy - kz,\\
{q^{\left( {\rm{2}} \right)}} = t + ix - jy - kz,\\
{q^{\left( {\rm{3}} \right)}} = t - ix - jy + kz.

\end{array}\]

\noindent Then, the norm of $q,$  denoted by $\left\| q \right\|,$  is defined by

\[{\left\| q \right\|^4} = \left| {q\,{q^{\left( {\rm{1}} \right)}}\,{q^{\left( {\rm{2}} \right)}}\,{q^{\left( {\rm{3}} \right)}}} \right| = \,\,\left| {\left[
{{{\left( {t + y} \right)}^2} + {{\left( {x + z} \right)}^2}} \right]\,\,\left[ {{{\left( {t - y} \right)}^2} + {{\left( {x - z} \right)}^2}} \right]} \right| \ge
0.\]

\section{HAMILTON MATRIX of COMMUTATIVE QUATERIONS}

\noindent Let $q = t + ix + jy + kz$  be a commutative quaternion. We will define the linear map  $\sigma $ as  $\sigma :H \to H$ such that $\sigma \left( x \right) = qx.$
\noindent This operator is called Hamilton operator. Using our newly defined operator and the basis $\left\{ {1,i,j,k} \right\}$  of $H,$  we can write

\begin{equation}\label{3}
\begin{array}{l}
\sigma \left( 1 \right) = \left( {t + ix + jy + kz} \right)1 = t + ix + jy + kz\\
\sigma \left( i \right) = \left( {t + ix + jy + kz} \right)i =  - x + it - jz + ky\\
\sigma \left( j \right) = \left( {t + ix + jy + kz} \right)j = y + iz + jt + kx\\
\sigma \left( k \right) = \left( {t + ix + jy + kz} \right)k =  - z + iy - jx + kt
\end{array}
\end{equation}

\noindent Then, we find the following matrix representation

\[\sigma \left( q \right) = \left( {\begin{array}{*{20}{c}}
t&{ - x}&y&{ - z}\\
x&{\,\,\,t}&z&{\,\,\,y}\\
y&{ - z}&t&{ - x}\\
z&{\,\,y}&x&{\,\,\,t}
\end{array}} \right).\]

\noindent The matrix $\eta \left( q \right)$  is called Hamilton matrix. The multiplication of two commutative quaternions \\ $q = t + x{e_1} + y{e_2} + z{e_3}$ and ${q_1} =
{t_1} + {x_1}{e_1} + {y_1}{e_2} + {z_1}{e_3}$  is given by the help of  Hamilton matrix

\[q{q_1} = {q_1}q = \sigma \left( {q_1} \right)q = \sigma \left( q \right){q_1}.\]

\noindent Now we define the function

\[\begin{array}{l}
\sigma :H \to M\\
\,\,\,\,\,\,\,\,q \to \sigma \left( q \right)
\end{array}\]

\noindent where

\[{\rm{M = }}\left\{ {\left( {\begin{array}{*{20}{c}}
t&{ - x}&y&{ - z}\\
x&{\,\,\,t}&z&{\,\,y}\\
y&{ - z}&t&{ - x}\\
z&{\,\,\,y}&x&{\,\,t}
\end{array}} \right):\,\,\,t,x,y,z \in R} \right\}.\]

\noindent The function $\sigma $  is bijection and it comprises the following properties

\[\begin{array}{l}
\sigma \left( {q{q_1}} \right) = \sigma \left( q \right)\sigma \left( q_1 \right)\\
\sigma \left( {q + q_1} \right) = \sigma \left( q \right) + \sigma \left( q_1 \right).
\end{array}\]

\noindent Thus, $\sigma $ is linear isomorphism. Also we obtain

\[\begin{array}{l}
\sigma \left( 1 \right) = \left( {\begin{array}{*{20}{c}}
1&0&0&0\\
0&1&0&0\\
0&0&1&0\\
0&0&0&1
\end{array}} \right) = {I_4};\,\,\,\,\,\,\,\,\sigma \left( {{e_1}} \right) = \left( {\begin{array}{*{20}{c}}
0&1&0&0\\
{ - 1}&0&0&0\\
0&0&0&1\\
0&0&{ - 1}&0
\end{array}} \right) = {E_1};\\
\\
\sigma \left( {{e_2}} \right) = \left( {\begin{array}{*{20}{c}}
0&0&1&0\\
0&0&0&1\\
1&0&0&0\\
0&1&0&0
\end{array}} \right) = {E_2};\,\,\,\,\,\sigma \left( {{e_3}} \right) = \left( {\begin{array}{*{20}{c}}
0&0&0&1\\
0&0&{ - 1}&0\\
0&1&0&0\\
{ - 1}&0&0&0
\end{array}} \right) = {E_3};
\end{array}\]

\noindent So that

\[\begin{array}{l}
\,\,\,\,\,\,\,\,\,\,\,\,\,\,\,\,\,\,\,\,\,\,\,\,\,\,\,\,\,\,\,\,\,{E_1}{E_1} =  - {I_4},\,\,\,{E_2}{E_2} = {I_4},\,\,\,{E_3}{E_3} =  - {I_4}\\
{E_1}{E_2} = {E_2}{E_1} = {E_3},\,\,{E_2}{E_3} = {E_3}{E_2} = {E_1},\,\,{E_1}{E_3} = {E_3}{E_1} =  - {E_2}.
\end{array}\]

\noindent Therefore, ${E_1},{E_2},{E_3}$ the matrices  are identical to commutative quaternions bases $i,j,k$ .

\begin{theorem}
\noindent If $q$  and $q_1$  are commutative quaternions and $\sigma $  is the Hamilton operator as defined before, then we have the following properties

\begin{description}
\,\,\,\,  \item[i.] $\sigma \left( {qq_1} \right) = \sigma \left( q \right)\sigma \left( q_1 \right),$
\,\,\,\,  \item[ii.] $\sigma \left( {\sigma \left( q \right)q_1} \right) = \sigma \left( q \right)\sigma \left( q_1 \right),$
\,\,\,\,  \item[iii.] $q = q_1 \Leftrightarrow \sigma \left( q_1 \right) = \sigma \left( q \right),$
\,\,\,\,  \item[iv.] $\sigma \left( {q + q_1} \right) = \sigma \left( q \right) + \sigma \left( q_1 \right),$
\,\,\,\,  \item[v.] $\sigma \left( {\lambda q} \right) = \lambda \sigma \left( q \right),\,\,for\,\lambda  \in R,$
\,\,\,\,  \item[vi.] $Trace\left( {\sigma \left( q \right)} \right) = q + {q^{\rm{(1)}}} + {q^{\rm{(2)}}} + \,{q^{\rm{(3)}}} = 4t,$
 \,\,\,\, \item[vii.] ${\left\| q \right\|^4} = \left| {\det \left( {\sigma \left( q \right)} \right)} \right| = \,\left| {\left[ {{{\left( {t + y} \right)}^2}
 - \alpha {{\left( {x + z} \right)}^2}} \right]\,\,\left[ {{{\left( {t - y} \right)}^2} - \alpha {{\left( {x - z} \right)}^2}} \right]} \right|.$
\end{description}
\end{theorem}
\noindent \begin{proof}
iii, iv, v and vi can be easily shown. For now we will prove i, ii,  and vii.
\noindent \\

\noindent \textbf{i.} Let $q = t + ix + jy + kz$ and $\,q_1 = {t_1} + i{x_1} + j{y_1} + k{z_1}.$  Then  $qq_1 = {r_0} + i{r_1} + j{r_2} + k{r_3}$ where  ${r_0} = t{t_1}
- x{x_1} + y{y_1} - z{z_1},\,\,\\ {r_1} = x{t_1} + t{x_1} + z{y_1} + y{z_1},\,\,{r_2} = t{y_1} + y{t_1} - x{z_1} - z{x_1}$ and  $\,{r_3} = z{t_1} + t{z_1} + x{y_1}
+ y{x_1}.$  Thus, we obtain
\noindent \\
\[\begin{array}{l}
\sigma \left( {qq_1} \right) = \left( {\begin{array}{*{20}{c}}
{{r_0}}&{ - {r_1}}&{{r_2}}&{ - {r_3}}\\
{{r_1}}&{{r_0}}&{{r_3}}&{{r_2}}\\
{{r_2}}&{ - {r_3}}&{{r_0}}&{ - {r_1}}\\
{{r_3}}&{{r_2}}&{{r_1}}&{{r_0}}
\end{array}} \right) = \left( {\begin{array}{*{20}{c}}
t&{ - x}&y&{ - z}\\
x&t&z&y\\
y&{ - z}&t&{ - x}\\
z&y&x&t
\end{array}} \right)\left( {\begin{array}{*{20}{c}}
{{t_1}}&{ - {x_1}}&{{y_1}}&{ - {z_1}}\\
{{x_1}}&{{t_1}}&{{z_1}}&{{y_1}}\\
{{y_1}}&{ - {z_1}}&{{t_1}}&{ - {x_1}}\\
{{z_1}}&{{y_1}}&{{x_1}}&{{t_1}}
\end{array}} \right)\\
\,\,\,\,\,\,\,\,\,\,\,\,\,\,\\
\,\,\,\,\,\,\,\,\,\,\,\,\,\,\,\,\, = \sigma \left( q \right)\sigma \left( q_1 \right).
\end{array}\]

\noindent \\

\noindent \textbf{ii.} \[\sigma \left( q \right)q_1 = \left( {\begin{array}{*{20}{c}}
t&{ - x}&y&{ - z}\\
x&t&z&y\\
y&{ - z}&t&{ - x}\\
z&y&x&t
\end{array}} \right)\left( {\begin{array}{*{20}{c}}
{{t_1}}\\
{{x_1}}\\
{{y_1}}\\
{{z_1}}
\end{array}} \right) = \left( {\begin{array}{*{20}{c}}
{{r_0}}\\
{{r_1}}\\
{{r_2}}\\
{{r_3}}
\end{array}} \right) = qq_1\]

\noindent Then, we have  $\sigma \left( {\sigma \left( q_1 \right)q} \right) = \sigma \left( q_1 \right)\sigma \left( q \right).$

\noindent \\
\noindent \textbf{vii.} \[\begin{array}{l}
\det \left( {\sigma \left( q_1 \right)} \right) = \det \left( {\begin{array}{*{20}{c}}
t&{ - x}&y&{ - z}\\
x&t&z&y\\
y&{ - z}&t&{ - x}\\
z&y&x&t
\end{array}} \right) = \det \left( {\begin{array}{*{20}{c}}
{t + y}&{ - x - z}&0&0\\
{x + z}&{t + y}&0&0\\
0&0&{t - y}&{ - x + z}\\
0&0&{x - z}&{t - y}
\end{array}} \right)\\
\\
\,\,\,\,\,\,\,\,\,\,\,\,\,\,\,\,\,\,\,\,\,\,\,\,\,\, = \left[ {{{\left( {t + y} \right)}^2} + {{\left( {x + z} \right)}^2}} \right]\left[ {{{\left( {t - y}
\right)}^2} + {{\left( {x - z} \right)}^2}} \right]
\end{array}\]

\noindent Then, we obtain ${\left\| q \right\|^4} = \left| {\det \left( {\sigma \left( p \right)} \right)} \right|.$

\end{proof}

\begin{theorem}
\noindent Every commutative quaternion can be represented by a  $2 \times 2$ complex matrix.

\end{theorem}

\noindent \begin{proof}
\noindent Let $q = t + ix + jy + kz \in H,$  then every commutative quaternion can be uniquely expressed as $q = {c_1} + j{c_2}\,,$  where ${c_1} = t + ix$  and ${c_2} = y +
iz$  are complex numbers. The linear map  $\varphi :H \to H$ is defined by ${\varphi _q}\left( p \right) = qp$  for all $p \in H.$  This map is bijective and

\[\begin{array}{l}
{\varphi _q}\left( 1 \right) = {c_1} + j{c_2}\\
{\varphi _q}\left( j \right) = {c_2} + j{c_1}
\end{array}\]

\noindent with this transformation,  commutative quaternions are defined as subset of the matrix ring ${M_2}\left( C \right),$  the set of $2 \times 2$  complex matrices:
\[N = \left\{ {\left( {\begin{array}{*{20}{c}}
{{c_1}}&{{c_2}}\\
{{c_2}}&{{c_1}}
\end{array}} \right):\,\,{c_1},{c_2} \in C} \right\}.\]

\noindent Then, $H$  and  $N$ are essentially the same.

\noindent Note that

\[\begin{array}{l}
\,\,\,\,\,\,\,\,\,\,\,\,\,\,\,\,\,\,\,\varphi :H \to N\\
\,\,\,q = {c_1} + j{c_2} \to \varphi \left( q \right) = \left( {\begin{array}{*{20}{c}}
{{c_1}}&{{c_2}}\\
{{c_2}}&{{c_1}}
\end{array}} \right)
\end{array}\]

\noindent is bijective and preserves the operations. Furthermore ${\left\| q \right\|^4} = \left| {\det \varphi \left( q \right)} \right|.$

\end{proof}

\section{COMMUTATIVE QUATERNION MATRICES}

\noindent The set ${M_{m \times n}}\left( H \right)$  denotes all $m \times n$  type matrices with commutative quaternion entries. If $m = n,$ then the set of commutative
quaternion matrices is denoted by ${M_n}\left( H \right).$  The ordinary matrix addition and multiplication are defined. The scalar multiplication is defined as, for $A = \left( {{a_{ij}}} \right) \in {M_{m \times n}}\left( H \right),\,\,q \in H,$

\[qA = Aq = \left( {q{a_{ij}}} \right).\]

\noindent It is easy to see that for $A \in {M_{m \times n}}\left( H \right)\,,\,\,\,B \in {M_{n \times k}}\left( H \right)\,{\rm{ and }}\,\,q,{q_1} \in H$

\[\begin{array}{l}
\left( {qA} \right)B = q\left( {AB} \right),\\
\left( {Aq} \right)B = A\left( {qB} \right),\\
\left( {q{q_1}} \right)A = q\left( {{q_1}A} \right).
\end{array}\]

\noindent Moreover, ${M_{m \times n}}\left( H \right)$  is free module over the ring $H.$
For $A = \left( {{a_{ij}}} \right) \in {M_{m \times n}}\left( H \right),$  there exists three kinds of its conjugates, called principal conjugates:

\[{A^{\left( {\rm{1}} \right)}} = \left( {a_{ij}^{\left( 1 \right)}} \right) \in {M_{m \times n}}\left( H \right),\,\,{A^{\left( 2 \right)}} = \left(
{a_{ij}^{\left( 2 \right)}} \right) \in {M_{m \times n}}\left( H \right)\,\,\,{\rm{and  }}\,\,\,{A^{\left( 3 \right)}} = \left( {a_{ij}^{\left( 3 \right)}} \right)
\in {M_{m \times n}}\left( H \right).\]

\noindent ${A^T} = \left( {{a_{ji}}} \right) \in {M_{m \times n}}\left( H \right)$  is the transpose of  $A\,;\,$ ${A^{{\dag _i}}} = {\left( {{A^{\left( {\rm{i}} \right)}}}
\right)^T} \in {M_{m \times n}}\left( H \right)$ is the ${{\mathop{\rm i}\nolimits} ^{th}}$ conjugate transpose of $A.$

\noindent A square matrix $A \in {M_n}\left( H \right)$  is said to be normal matrix by ${i^{th}}$   conjugate, if $A{A^{{\dag _i}}} = {A^{{\dag _i}}}A;$ Hermitian matrix  by
${i^{th}}$  conjugate,  if $A = {A^{{\dag _i}}};$  unitary matrix by ${i^{th}}$   conjugate ,  if $A{A^{{\dag _i}}} = I$  where $I$ is the identity matrix and
invertible matrix, if  $AB = BA = I$  for some $B \in {M_n}\left( H \right).$

\begin{theorem}

Let $A \in {M_{m \times n}}\left( H \right),\,\,B \in {M_{n \times p}}\left( H \right).$  Then  the followings are satisfied:

\begin{description}
  \item[I.] ${\left( {{A^{\left( i \right)}}} \right)^T} = {\left( {{A^T}} \right)^{\left( i \right)}}\,;$
  \item[II.] ${\left( {AB} \right)^{{\dag _i}}} = {B^{{\dag _i}}}{A^{{\dag _i}}}\,;$
  \item[III.] ${\left( {AB} \right)^T} = {B^T}{A^T}\,;$
  \item[IV.] ${\left( {AB} \right)^{\left( i \right)}} = {A^{\left( i \right)}}{B^{\left( i \right)}}\,;$
  \item[V.] ${\left( {AB} \right)^{ - 1}} = {B^{ - 1}}{A^{ - 1}}\,$ if $A$    and $B$  is invertible;
  \item[VI.] ${\left( {{A^{{\dag _i}}}} \right)^{ - 1}} = {\left( {{A^{ - 1}}} \right)^{{\dag _i}}}$ if  $A$ is invertible;
  \item[VII.] ${\left( {{A^{\rm{(i)}}}} \right)^{ - 1}} = {\left( {{A^{ - 1}}} \right)^{\rm{(i)}}};$
  \item[VIII.] ${\left( {{A^{\left( i \right)}}} \right)^{\left( j \right)}} = \left\{ {\begin{array}{*{20}{c}}
{{A^{\left( k \right)}}}&{}&{{\rm{i}} \ne {\rm{j}} \ne {\rm{k}}}\\
{}&{}&{}\\
A&{}&{{\rm{i}} = {\rm{j}}}
\end{array}} \right.\,\,\,\,\,\,\,\,\,\,\,{\rm{for  }}\,\,\,i,j,k \in \left\{ {1,2,3} \right\};$
\end{description}

\noindent where $A^{\left( i \right)}$  is the  ${i^{th}}$ conjugate and $i = 1,2,3.$

\end{theorem}

\noindent \begin{proof}

\noindent I, III, IV, V, VI, VII and VIII can be easily shown. Now we will prove II.

\noindent \textbf{II.} Let $A = {A_1} + j{A_2}$  and $B = {B_1} + j{B_2},$  where ${A_1},{A_2},{B_1}$  and ${B_2}$  are complex matrices. Then

\[\begin{array}{l}
{\left( {AB} \right)^{{\dag _1}}} = {\left[ {\left( {{A_1} + j{A_2}} \right)\left( {{B_1} + j{B_2}} \right)} \right]^{{\dag _1}}}\\
\,\,\,\,\,\,\,\,\,\,\,\,\,\,\, = {\left[ {\overline {\left( {{A_1}{B_1}} \right)}  + \overline {\left( {{A_2}{B_2}} \right)}  + j\left[ {\overline {\left(
{{A_1}{B_2}} \right)}  + \overline {\left( {{A_2}{B_1}} \right)} } \right]} \right]^T}\\
\,\,\,\,\,\,\,\,\,\,\,\,\,\,\, = {\left[ {{{\overline A }_1}\overline B _1^{} + \overline A _2^{}\overline B _2^{} + j\left( {\overline A _1^{}\overline B _2^{} +
\overline A _2^{}\overline B _1^{}} \right)} \right]^T}\\
\,\,\,\,\,\,\,\,\,\,\,\,\,\,\, = {\left( {\overline B _1^{}} \right)^T}{\left( {\overline A _1^{}} \right)^T} + {\left( {\overline B _2^{}} \right)^T}{\left(
{\overline A _2^{}} \right)^T} + j\left[ {{{\left( {\overline B _2^{}} \right)}^T}{{\left( {\overline A _1^{}} \right)}^T} + {{\left( {\overline B _1^{}}
\right)}^T}{{\left( {\overline A _2^{}} \right)}^T}} \right]
\end{array}\]

\noindent On the other hand, for

\[\begin{array}{l}
{B^{{\dag _1}}} = {\left( {\overline B _1^{}} \right)^T} + j{\left( {\overline B _2^{}} \right)^T}\,\,\,\,\,{\rm{and}}\,\,\,\,\,\,{A^{{\dag _1}}} = {\left(
{\overline A _1^{}} \right)^T} + j{\left( {\overline A _2^{}} \right)^T}\\
\\
{B^{{\dag _1}}}{A^{{\dag _1}}} = {\left( {\overline B _1^{}} \right)^T}{\left( {\overline A _1^{}} \right)^T} + {\left( {\overline B _2^{}} \right)^T}{\left(
{\overline A _2^{}} \right)^T} + j\left[ {{{\left( {\overline B _2^{}} \right)}^T}{{\left( {\overline A _1^{}} \right)}^T} + {{\left( {\overline B _1^{}}
\right)}^T}{{\left( {\overline A _2^{}} \right)}^T}} \right].
\end{array}\]

\noindent Thus, we obtain ${\left( {AB} \right)^{{\dag _1}}} = {B^{{\dag _1}}}{A^{{\dag _1}}}$  .

\[\begin{array}{l}
{\left( {AB} \right)^{{\dag _2}}} = {\left[ {\left( {{A_1} + j{A_2}} \right)\left( {{B_1} + j{B_2}} \right)} \right]^{{\dag _2}}}\\
\,\,\,\,\,\,\,\,\,\,\,\,\,\,\, = {\left[ {\left( {{A_1}{B_1}} \right) + \left( {{A_2}{B_2}} \right) - j\left[ {\left( {{A_1}{B_2}} \right) + \left( {{A_2}{B_1}}
\right)} \right]} \right]^T}\\
\,\,\,\,\,\,\,\,\,\,\,\,\,\,\, = \left[ {B_1^TA_1^T + B_2^TA_2^T - j\left( {B_2^TA_1^T + B_1^TA_2^T} \right)} \right]
\end{array}\]

\noindent On the other hand , for

\[\begin{array}{l}
{B^{{\dag _2}}} = B_1^T - jB_2^T\,\,\,\,{\rm{and}}\,\,\,\,\,\,{A^{{\dag _1}}} = A_1^T - jA_2^T\\
\\
{B^{{\dag _2}}}{A^{{\dag _2}}} = \left[ {B_1^TA_1^T + B_2^TA_2^T - j\left( {B_2^TA_1^T + B_1^TA_2^T} \right)} \right].
\end{array}\]

\noindent Thus, we have ${\left( {AB} \right)^{{\dag _2}}} = {B^{{\dag _2}}}{A^{{\dag _2}}}.$  The same way, we have ${\left( {AB} \right)^{{\dag _3}}} = {B^{{\dag
_3}}}{A^{{\dag _3}}}.$

\end{proof}

\begin{theorem}

Let $A,B \in {M_n}\left( H \right),\,$ if $AB = I$  then $BA = I.$

\end{theorem}

\noindent \begin{proof}

\noindent Fist note that the proposition is true for complex matrices.  Let $A = {A_1} + j{A_2}$  and $B = {B_1} + j{B_2},$  where ${A_1},{A_2},{B_1}$ and ${B_2}$  are $n
\times n$  complex matrices and ${j^2} = 1$ . Then,

\[AB = \left( {{A_1}{B_1} + {A_2}{B_2}} \right) + j\left( {{A_1}{B_2} + {A_2}{B_1}} \right) = I.\]

\noindent From this we have

\[\left( {\begin{array}{*{20}{c}}
{{A_1}}&{{A_2}}
\end{array}} \right)\left( {\begin{array}{*{20}{c}}
{{B_1}}&{{B_2}}\\
{{B_2}}&{{B_1}}
\end{array}} \right) = \left( {\begin{array}{*{20}{c}}
I&0
\end{array}} \right)\]

\noindent so

\[\left( {\begin{array}{*{20}{c}}
{{A_1}}&{{A_2}}\\
{{A_2}}&{{A_1}}
\end{array}} \right)\left( {\begin{array}{*{20}{c}}
{{B_1}}&{{B_2}}\\
{{B_2}}&{{B_1}}
\end{array}} \right) = \left( {\begin{array}{*{20}{c}}
I&0\\
0&I
\end{array}} \right).\]

\noindent Since $\left( {\begin{array}{*{20}{c}}
{{A_1}}&{{A_2}}\\
{{A_2}}&{{A_1}}
\end{array}} \right)\,\,\,\,{\rm{and}}\,\,\,\left( {\begin{array}{*{20}{c}}
{{B_1}}&{{B_2}}\\
{{B_2}}&{{B_1}}
\end{array}} \right)$  are $2n \times 2n$  complex  matrices, we  get

\[\left( {\begin{array}{*{20}{c}}
{{B_1}}&{{B_2}}\\
{{B_2}}&{{B_1}}
\end{array}} \right)\left( {\begin{array}{*{20}{c}}
{{A_1}}&{{A_2}}\\
{{A_2}}&{{A_1}}
\end{array}} \right) = \left( {\begin{array}{*{20}{c}}
I&0\\
0&I
\end{array}} \right).\]

\noindent Then we can write  ${B_1}{A_1} + {B_2}{A_2} = I,\,\,{B_1}{A_2} + {B_2}{A_1} = 0.$   So $\left( {{B_1}{A_1} + {B_2}{A_2}} \right) + j{B_1}{A_2} + {B_2}{A_1} =
I.$ Thus we obtain  $BA = I.$

\end{proof}

\begin{definition}
Let $A \in {M_n}\left( H \right)$  and $\lambda  \in H.$  If $\lambda $  satisfies the equation of  $Ax = \lambda x$ then $\lambda $ is called the eigenvalue of
$A.$ The set of the eigenvalues is defined as

\[\xi \left( A \right) = \left\{ {\lambda  \in H:\,\,Ax = \lambda x\,\,{\rm{for some }}x \ne 0} \right\}.\]

\noindent This set is called the spectrum of $A.$

\end{definition}

\noindent We shall always write every commutative quaternion $q = t + ix + jy + kz$  in the form $q = {c_1} + j{c_2}$  where ${c_1} = t + ix$  and  ${c_2} = y + iz$ are
special commutative quaternions. The  special commutative quaternions, i.e. quaternions of the form $\alpha  + i\beta \,\,\,\left( {\alpha ,\beta  \in R} \right),$
form a commutative field $C$ which is algebraically closed. Since $C$  is evidently isomorphic to a field of ordinary complex numbers, we shall call the special
commutative quaternions in $C$  simply as complex numbers.

\noindent Let $A$  be an $n \times n$  matrix with coefficients in $H.$  We shall show that there exists a non-zero  column vector  $\left( {n \times 1\,\,{\rm{matrix}}}
\right)$ $x$   with components in $H$  such that

\begin{equation}\label{4}
Ax = \lambda x\,,\,\,\,\lambda  \in C
\end{equation}

\noindent The multiplier $\lambda ,$  which we restrict to be a complex number ( special commutative quaternion), will be called an eigenvalue of the matrix $A.$

\noindent Decompose $A = {A_1} + j{A_2},\,\,x = {x_1} + j{x_2},$  where ${A_1},{A_2}$  are  $n \times n$ complex matrices and ${x_1},{x_2}$  are column vectors with
components in $C.$  Then (\ref{4}) maybe written

\[\left( {{A_1} + j{A_2}} \right)\left( {{x _1} + j{x _2}} \right) = \lambda {x_1} + j\lambda {x_2}\]

\noindent hence

\[{A_1}{x_1} + {A_2}{x_2} = \lambda {x_1}\,\,,\,\,\,\,{A_1}{x_2} + {A_2}{x_1} = \lambda {x_2}\,\]

\noindent which are equations in the commutative complex field $C$  and may be written in the matrix form

\begin{equation}\label{5}
\left( {\begin{array}{*{20}{c}}
{{A_1}}&{{A_2}}\\
{{A_2}}&{{A_1}}
\end{array}} \right)\left( {\begin{array}{*{20}{c}}
{{x_1}}\\
{{x_2}}
\end{array}} \right) = \lambda \left( {\begin{array}{*{20}{c}}
{{x_1}}\\
{{x_2}}
\end{array}} \right).
\end{equation}

\noindent Since $C$  is algebraically closed, (\ref{5}) implies the existence of $2n$  values of $\lambda$  to each of which correspond a non-zero vector $\left(
{\begin{array}{*{20}{c}}
{{x_1}}\\
{{x_2}}
\end{array}} \right)$    satisfying (\ref{5}). and therefore a non-zero vector $\lambda $ satisfying (\ref{4}). Thus,

\begin{theorem}
\noindent An  $n \times n$ matrix $A$  with commutative quaternion coefficients has $2n$  eigenvalues which are complex numbers.

\end{theorem}

\section{THE COMPLEX ADJOINT MATRICES OF THE COMMUTATIVE QUATERNION MATRICES}

\noindent In this section, we will define the complex adjoint matrix of a commutative quaternion matrix. After that we will give some relations between commutative quaternion
matrices and their complex adjoint matrices.

\begin{definition}
Let $A = {A_1} + j{A_2} \in {M_n}\left( H \right)$  where  ${A_1}$ and ${A_2}$  are complex matrices. We define the $2n \times 2n$  complex  matrices  $\eta \left(
A \right)$ as

\[\left( {\begin{array}{*{20}{c}}
{{A_1}}&{{A_2}}\\
{{A_2}}&{{A_1}}
\end{array}} \right).\]

\noindent This matrix $\eta \left( A \right)$  is called the complex adjoint (or adjoint) matrix of the commutative quaternion matrix $A$ .

\end{definition}

\begin{definition}
Let $A \in {M_n}\left( H \right)$  and $\eta \left( A \right)$  be the adjoint matrix of $A.$  We define the $q-$ determinant of $A$  by \\ ${\det _q}\left( A \right)
= \det \left( {\eta \left( A \right)} \right).$  Here $\det \left( {\eta \left( A \right)} \right)$  is usual determinant of $\eta \left( A \right).$

\end{definition}

\noindent For example Let $A = \left( {\begin{array}{*{20}{c}}
i&{j + k}\\
1&k
\end{array}} \right)$  be a commutative quaternion matrix. Then the adjoint matrix of $A$  is

\[\eta \left( A \right) = \left( {\begin{array}{*{20}{c}}
i&0&0&{1 + i}\\
1&0&0&i\\
0&{1 + i}&i&0\\
0&i&1&0
\end{array}} \right).\]

\noindent ${\det _q}\left( A \right) = \det \left( {\eta \left( A \right)} \right) =  - 3 - 4i.$  Note that $\eta \left( {{A^{\left( 1 \right)}}} \right) = \overline {\eta
\left( A \right)} .$  But $\eta \left( {{A^{\left( 2 \right)}}} \right) \ne \overline {\eta \left( A \right)} $   and  $\eta \left( {{A^{\left( 3 \right)}}} \right)
\ne \overline {\eta \left( A \right)} $ in general, for $A \in {M_n}\left( H \right).$

\begin{theorem}

Let $A,B \in {M_n}\left( H \right).$  Then the followings hold

\begin{description}
  \item[I.] $\eta \left( {{I_n}} \right) = {I_{2n}}\,;$
  \item[II.] $\eta \left( {A + B} \right) = \eta \left( A \right) + \eta \left( B \right);$
  \item[III.] $\eta \left( {AB} \right) = \eta \left( A \right)\eta \left( B \right);$
  \item[IV.] $\eta \left( {{A^{ - 1}}} \right) = {\left( {\eta \left( A \right)} \right)^{ - 1}}$ if ${A^{ - 1}}$ exist;
  \item[V.] $\eta \left( {{A^{{\dag _1}}}} \right) = {\left( {\eta \left( A \right)} \right)^\dag },$ but $\eta \left( {{A^{{\dag _2}}}} \right) \ne {\left(
      {\eta \left( A \right)} \right)^\dag },\eta \left( {{A^{{\dag _3}}}} \right) \ne {\left( {\eta \left( A \right)} \right)^\dag }$ in general,\, ${\left(
      {\eta \left( A \right)} \right)^\dag }$ is the conjugate transpose  of $\eta \left( A \right);$
  \item[VI.] ${\det _q}\left( {AB} \right) = {\det _q}\left( A \right){\det _q}\left( B \right),$ consequently ${\det _q}\left( {{A^{ - 1}}} \right) = {\left(
      {{{\det }_q}\left( A \right)} \right)^{ - 1}},$ if ${A^{ - 1}}$  is exists;
\end{description}

\end{theorem}

\noindent \begin{proof}
\noindent I, II and IV can be easily shown. For now we will prove III , V and VI.

\noindent \textbf{III.} Let $A,B \in {M_n}\left( H \right),\,$ where $A = {A_1} + j{A_2}$  and $B = {B_1} + j{B_2}$  for the complex matrices  ${A_1},{A_2},{B_1}$ and
      ${B_2}.$ Then the adjoint matrices of $A$  and $B$  are

\[\eta \left( A \right) = \left( {\begin{array}{*{20}{c}}
{{A_1}}&{{A_2}}\\
{{A_{\,2}}}&{{A_1}}
\end{array}} \right),\,\,\,\,\,\,\eta \left( B \right) = \left( {\begin{array}{*{20}{c}}
{{B_1}}&{{B_2}}\\
{{B_{\,2}}}&{{B_1}}
\end{array}} \right).\]

\noindent If we calculate $AB$  for $A = {A_1} + j{A_2}$  and $B = {B_1} + j{B_2},$ then we have

\[AB = \left( {{A_1}{B_1} + {A_2}{B_2}} \right) + j\left( {{A_1}{B_2} + {A_2}{B_2}} \right).\]

\noindent So far the adjoint matrix of $AB$ we can write

\[\eta \left( A \right) = \left( {\begin{array}{*{20}{c}}
{{A_1}{B_1} + {A_2}{B_2}}&{{A_1}{B_2} + {A_2}{B_2}}\\
{{A_1}{B_2} + {A_2}{B_2}}&{{A_1}{B_1} + {A_2}{B_2}}
\end{array}} \right)\]
\noindent On the other hand,

\[\eta \left( A \right)\eta \left( B \right) = \left( {\begin{array}{*{20}{c}}
{{A_1}}&{{A_2}}\\
{{A_{\,2}}}&{{A_1}}
\end{array}} \right)\left( {\begin{array}{*{20}{c}}
{{B_1}}&{{B_2}}\\
{{B_{\,2}}}&{{B_1}}
\end{array}} \right) = \left( {\begin{array}{*{20}{c}}
{{A_1}{B_1} + {A_2}{B_2}}&{{A_1}{B_2} + {A_2}{B_2}}\\
{{A_1}{B_2} + {A_2}{B_2}}&{{A_1}{B_1} + {A_2}{B_2}}
\end{array}} \right).\]

\noindent Thus we have  $\eta \left( {AB} \right) = \eta \left( A \right)\eta \left( B \right).$

\noindent \textbf{V.} Let $A = {A_1} + j{A_2} \in {M_n}\left( H \right).$  Since ${A^{{\dag _1}}} = {\left( {\overline {{A_1}} } \right)^T} + j{\left( {\overline {{A_2}}
      } \right)^T},$   adjoint matrix of ${A^{{\dag _1}}}$  is

      \[\eta \left( {{A^{{\dag _1}}}} \right) = \left( {\begin{array}{*{20}{c}}
{{{\left( {\overline {{A_1}} } \right)}^T}}&{{{\left( {\overline {{A_2}} } \right)}^T}}\\
{{{\left( {\overline {{A_2}} } \right)}^T}}&{{{\left( {\overline {{A_2}} } \right)}^T}}
\end{array}} \right)\]

\noindent On the other hand

\[\eta \left( A \right) = \left( {\begin{array}{*{20}{c}}
{{A_1}}&{{A_2}}\\
{{A_2}}&{{A_1}}
\end{array}} \right)\,\,{\rm{   and    }}\,\,{\left( {\eta \left( A \right)} \right)^\dag } = \left( {\begin{array}{*{20}{c}}
{{{\left( {\overline {{A_1}} } \right)}^T}}&{{{\left( {\overline {{A_2}} } \right)}^T}}\\
{{{\left( {\overline {{A_2}} } \right)}^T}}&{{{\left( {\overline {{A_2}} } \right)}^T}}
\end{array}} \right).\]

\noindent Thus we obtain  $\eta \left( {{A^{{\dag _1}}}} \right) = {\left( {\eta \left( A \right)} \right)^\dag }.$

\noindent Moreover

    \[\begin{array}{l}
\eta \left( {{A^{{\dag _2}}}} \right) = \left( {\begin{array}{*{20}{c}}
{A_1^T}&{ - A_2^T}\\
{ - A_2^T}&{A_1^T}
\end{array}} \right) \ne \left( {\begin{array}{*{20}{c}}
{{{\left( {\overline {{A_1}} } \right)}^T}}&{{{\left( {\overline {{A_2}} } \right)}^T}}\\
{{{\left( {\overline {{A_2}} } \right)}^T}}&{{{\left( {\overline {{A_2}} } \right)}^T}}
\end{array}} \right)\,\,{\rm{   in\,\, general}}\\
\\
\eta \left( {{A^{{\dag _3}}}} \right) = \left( {\begin{array}{*{20}{c}}
{{{\left( {\overline {{A_1}} } \right)}^T}}&{ - {{\left( {\overline {{A_2}} } \right)}^T}}\\
{ - {{\left( {\overline {{A_2}} } \right)}^T}}&{{{\left( {\overline {{A_1}} } \right)}^T}}
\end{array}} \right) \ne \left( {\begin{array}{*{20}{c}}
{{{\left( {\overline {{A_1}} } \right)}^T}}&{{{\left( {\overline {{A_2}} } \right)}^T}}\\
{{{\left( {\overline {{A_2}} } \right)}^T}}&{{{\left( {\overline {{A_2}} } \right)}^T}}
\end{array}} \right)\,\,{\rm{    in\,\, general}}
\end{array} \]

\noindent \textbf{VI.} From Theorem (9-III) we obtain $\eta \left( {AB} \right) = \eta \left( A \right)\eta \left( B \right).$  And so

      \[\begin{array}{l}
{\det _q}\left( {AB} \right) = \det \left( {\eta \left( {AB} \right)} \right) = \det \left( {\eta \left( A \right)\eta \left( B \right)} \right)\\
\,\,\,\,\,\,\,\,\,\,\,\,\,\,\,\,\,\,\,\,\,\, = \det \left( {\eta \left( A \right)} \right)\det \left( {\eta \left( B \right)} \right) = {\det _q}\left( A
\right){\det _q}\left( B \right).
\end{array}\]

\noindent Here, if we let $B = {A^{ - 1}}$  then we can find easily  ${\det _q}\left( {{A^{ - 1}}} \right) = {\left( {{{\det }_q}\left( A \right)} \right)^{ - 1}}.$

\end{proof}

\begin{theorem}

\noindent Let $A \in {M_n}\left( H \right).$  The following are equivalent:

\begin{description}
  \item[(I.)] $A$ is invertible,
  \item[(II.)] $Ax = 0$ has a unique solution,
  \item[(III.)] $\det \left( {\eta \left( A \right)} \right) \ne 0$ i.e $\eta \left( A \right)$ is an invertible,
  \item[(IV.)] $A$ has no zero eigenvalue. More precisely, if $Ax = \lambda x$  for some commutative quaternion $\lambda $ and some commutative vector $x \ne
      0,$$x \ne 0,$  then $\lambda  \ne 0.$
\end{description}

\end{theorem}

\noindent \begin{proof}
$\left( {\textbf{I.}} \right) \Rightarrow \left( {\textbf{II.}} \right)$ this is trivial.

\noindent $\left( {\textbf{II.}} \right) \Rightarrow \left( {\textbf{III.}} \right)$ Let $A = {A_1} + j{A_2},\,\,x = {x_1} + j{x_2}$  where  ${A_1},{A_2}$ are complex  matrices and
${x_1},{x_2}$  are complex column vectors. Then

\[Ax = \left( {{A_1}{x_1} + {A_2}{x_2}} \right) + j\left( {{A_1}{x_2} + {A_2}{x_1}} \right).\]

\noindent From $Ax = 0$  we can write

\[\left( {{A_1}{x_1} + {A_2}{x_2}} \right) = 0\,\,\,\,\,\,{\rm{and}}\,\,\,\,\,\left( {{A_1}{x_2} + {A_2}{x_1}} \right) = 0.\]

\noindent So we get that

\[Ax = 0\,\,\,\,\,{\rm{if}}\,\,{\rm{and}}\,\,{\rm{only}}\,\,{\rm{if}}\,\,\,\,\,\,\left( {\begin{array}{*{20}{c}}
{{A_1}}&{{A_2}}\\
{{A_2}}&{{A_1}}
\end{array}} \right)\left( {\begin{array}{*{20}{c}}
{{x_1}}\\
{{x_2}}
\end{array}} \right) = 0.\]

\noindent That is $\eta \left( A \right){\left( {\begin{array}{*{20}{c}}
{{x_1}}&{{x_2}}
\end{array}} \right)^T} = 0.$  Since $Ax = 0$  has an unique solution,$\eta \left( A \right){\left( {\begin{array}{*{20}{c}}
{{x_1}}&{{x_2}}
\end{array}} \right)^T} = 0$   has  an unique solution. Thus, since $\eta \left( A \right)$  is complex  matrix, $\eta \left( A \right)$  is invertible.

\noindent $\left( {\textbf{II.}} \right) \Rightarrow \left( {\textbf{IV.}} \right)$ Let  $Ax = 0$ has a unique solution $0$  for $A \in {M_n}\left( H \right).$  Suppose that $A$  has no a zero
eigenvalue. Then for some commutative quaternion vector $x = 0,$ the equation $Ax = \lambda x$  has zero eigenvalue. Thus   and so by our assumption $Ax = 0$  this
is a contradiction. Now suppose that   has no zero eigenvalue. If we have $Ax = 0 = \lambda x,$  then by our assumption $x = 0.$

\noindent $\left( {\textbf{III.}} \right) \Rightarrow \left( {\textbf{I.}} \right)$ If $\eta \left( A \right)$  is invertible, then for  $A = {A_1} + j{A_2}$ there exist a complex matrix
$\left( {\begin{array}{*{20}{c}}
{{B_1}}&{{B_2}}\\
{{B_2}}&{{B_1}}
\end{array}} \right)$  such that

\[\left( {\begin{array}{*{20}{c}}
{{B_1}}&{{B_2}}\\
{{B_2}}&{{B_1}}
\end{array}} \right)\left( {\begin{array}{*{20}{c}}
{{A_1}}&{{A_2}}\\
{{A_2}}&{{A_1}}
\end{array}} \right) = \left( {\begin{array}{*{20}{c}}
I&0\\
0&I
\end{array}} \right).\]

\noindent Here we obtain  ${B_1}{A_1} + {B_2}{A_2} = I\,\,{\rm{and}}\,\,\,{B_1}{A_2} + {B_2}{A_1} = 0.$  Using this equation we can write

\[\left( {{B_1}{A_1} + {B_2}{A_2}} \right) + j\left( {{B_1}{A_2} + {B_2}{A_1}} \right) = I.\]

\noindent That is $BA = I$  for  $B = {B_1} + j{B_2}.$ So $A$  is invertible commutative quaternion matrix by Theorem (4).

\end{proof}



%% References



\end{document}